\documentclass[11pt]{amsart}
\usepackage[utf8]{inputenc}
\usepackage[T1]{fontenc}
\usepackage{amsmath,amssymb,amsthm}
\usepackage{booktabs}
\usepackage{enumerate}
\usepackage[margin=2.7cm]{geometry}
\usepackage[colorlinks=true,linkcolor=blue,citecolor=blue,urlcolor=blue]{hyperref}

\newtheorem{thm}{Theorem}[section]
\newtheorem{prop}[thm]{Proposition}
\newtheorem{lem}[thm]{Lemma}
\newtheorem{cor}[thm]{Corollary}
\newtheorem{conj}[thm]{Conjecture}

\theoremstyle{definition}
\newtheorem{defn}[thm]{Definition}
\newtheorem{ex}[thm]{Example}
\newtheorem{rk}[thm]{Remark}

\newcommand{\C}{\mathbb{C}}
\newcommand{\PP}{\mathbb{P}}
\newcommand{\mdr}{\operatorname{mdr}}
\newcommand{\AR}{\operatorname{AR}}
\newcommand{\Sing}{\operatorname{Sing}}
\newcommand{\Syz}{\operatorname{Syz}}

\begin{document}

\title[Plane curves attaining the Dimca--Saito--Sernesi bound]{Plane curves attaining the Dimca--Saito--Sernesi bound and Jacobian syzygies of pencil-type curves}

\author{Gabriel Sticlaru}
\address{Faculty of Mathematics and Informatics, Ovidius University, Bd.\ Mamaia 124, 900527 Constan\c{t}a, Romania}
\email{gabriel.sticlaru@gmail.com}


\begin{abstract}
For a reduced plane curve $C:f=0$ of degree $d$ with only weighted homogeneous singularities, Dimca, Saito and Sernesi proved that the minimal degree of a Jacobian relation satisfies $\mdr(f)\geq \alpha_C\, d-2$, where $\alpha_C$ is the minimal Arnold exponent of the singularities of $C$. We call $C$ \emph{DSS-extremal} (Dimca--Saito--Sernesi extremal) when equality holds. We show that curves of pencil type $f=\Phi(g,h)$, in particular curves of $(p,q)$-torus type $g^p+h^q$, are DSS-extremal, the minimal syzygy being $\nabla g\times\nabla h$, and that this syzygy never violates the bound, whatever the singularities. Under a natural transversality assumption we determine the minimal free resolution of the Milnor algebra of such a curve: it depends only on $d$, $\deg g$ and $\deg h$, and the module of Jacobian syzygies is generated by $\nabla g\times \nabla h$ and the three Koszul syzygies. For irreducible sextics with six cusps we prove that $3\le \mdr(f)\le 4$, that $\mdr(f)=3$ for curves of torus type, and that $\mdr(f)=4$ for the general curve in the component of non-torus sextics containing Oka's curve. Hence Zariski's classical pair is a Ziegler pair detected by $\mdr$, and the reason is the extremality of the torus member. We conjecture that every DSS-extremal curve is of pencil type and support this by the nodal case and by computations. Finally, we conjecture that the Dimca--Saito--Sernesi bound holds for all reduced plane curves, with $\alpha_C$ the minimal Arnold exponent, and we prove it for curves with $\mdr(f)\le1$. These include an infinite family of free curves, found by du Plessis and Wall, whose unique singularity is not weighted homogeneous.
\end{abstract}

\maketitle

\section{Introduction}

Let $S=\C[x,y,z]$ be the graded polynomial ring and let $C:f=0$ be a reduced curve of degree $d\ge 3$ in the complex projective plane $\PP^2$. Denote by $J_f=(f_x,f_y,f_z)$ the Jacobian ideal, by $M(f)=S/J_f$ the Milnor (Jacobian) algebra and by
\[
\AR(f)=\{(a,b,c)\in S^3 \ : \ af_x+bf_y+cf_z=0\}
\]
the graded module of Jacobian syzygies. The \emph{minimal degree of a Jacobian relation} $\mdr(f)$ is the smallest $r$ with $\AR(f)_r\neq 0$. This invariant controls the freeness of $C$ through the bounds of du Plessis and Wall \cite{dPW}, and it is the basic invariant in the study of free and nearly free curves \cite{DS-RIMS,DS-3syz,Cuntz-Pokora}, of Terao's conjecture and of Ziegler pairs \cite{Ziegler,Cuntz-Pokora,Dimca-Pokora}.

When all the singularities of $C$ are weighted homogeneous, one has the lower bound
\begin{equation}\label{eq:DSS}
\mdr(f)\ \ge\ \alpha_C\, d-2,
\end{equation}
where $\alpha_C$ is the minimum of the Arnold exponents (log canonical thresholds) of the singular points of $C$. This inequality was proved by Dimca and Saito \cite[Corollary~5.5]{DSa}, and it was stated explicitly for plane curves, with this formulation, by Dimca and Sernesi \cite[Theorem~2.1]{DSe}; we call it the \emph{Dimca--Saito--Sernesi bound}. The bound \eqref{eq:DSS} is used systematically in recent work on free curves and curve arrangements. Dimca and Sernesi noted that it is sharp for the line arrangement $(x^2-y^2)(y^2-z^2)(z^2-x^2)=0$, for a nonic with nine $D_4$ points, and for Zariski's sextic with six cusps on a conic. To our knowledge, the curves for which equality holds have not been characterized.

\begin{defn}\label{def:extremal}
A reduced plane curve $C:f=0$ with only weighted homogeneous singularities is called \emph{Dimca--Saito--Sernesi extremal}, or \emph{DSS-extremal} for short, if $\mdr(f)=\alpha_C\, d-2$.
\end{defn}

This is a refinement in the same spirit as the passage from free to nearly free curves: among all curves with given singularities, the DSS-extremal ones have Jacobian syzygies in the lowest degree allowed by the local data. For curves with simple elliptic singularities, the free curves attaining \eqref{eq:DSS} were studied in \cite{M-curves} under the name of $M$-curves. We use the prefix DSS because the word \emph{extremal} has other meanings in algebraic geometry.

The examples above share a common feature. Each of them can be written as $f=\Phi(g,h)$ for two forms $g,h$ and a binary form $\Phi$ (Definition~\ref{def:pencil}). For such a curve the vector
\[
\rho=\nabla g\times \nabla h
\]
is a Jacobian syzygy of degree $\deg g+\deg h-2$, and the local equation at a base point $g=h=0$ is $\Phi$ itself. Our main results are the following.

\begin{enumerate}[(A)]
\item \emph{Pencil-type curves are DSS-extremal} (Theorem~\ref{thm:A}). If $f=\Phi(g,h)$, then $\alpha_C\le(\deg g+\deg h)/d$ for any base locus and any singularities (Proposition~\ref{prop:pencil-lct}), so the syzygy $\rho$ never violates \eqref{eq:DSS}. If the singularities are weighted homogeneous and equality holds, for instance when the minimal Arnold exponent of $C$ is attained at a transversal base point, then $\mdr(f)=\deg g+\deg h-2=\alpha_C d-2$. This covers the curves of $(p,q)$-torus type $g^p+h^q$ and unions of members of a pencil, for instance the Ceva and Hesse line arrangements (Example~\ref{ex:arr}).
\item \emph{The minimal resolution} (Theorem~\ref{thm:B}). If, in addition, $\rho$ does not vanish at the singular points of $C$, then $\AR(f)$ is minimally generated by $\rho$ and the three Koszul syzygies. There are exactly two second syzygies, and the graded Betti numbers of $M(f)$ depend only on $(d,\deg g,\deg h)$. Moreover $\tau(C)=(d-\deg g)(d-\deg h)$.
\item \emph{Sextics with six cusps} (Theorem~\ref{thm:C}). For such a sextic, $3\le\mdr(f)\le 4$. Torus-type sextics have $\mdr(f)=3$, while Oka's non-torus sextic, and the general curve of its component, have $\mdr(f)=4$. Hence Zariski's pair is a Ziegler pair distinguished by $\mdr$.
\item \emph{Beyond weighted homogeneity} (Section~\ref{sec:nonWH}). We conjecture that \eqref{eq:DSS} holds for every reduced curve (Conjecture~\ref{conj:nonWH}), and we prove it when $\mdr(f)\le1$ (Proposition~\ref{prop:mdr1}). Then the only curves with a singularity which is not weighted homogeneous are the curves $C_{k,\varepsilon}$ of du Plessis and Wall \cite{dPW99}, for which $\mdr(f)=1$ and $\alpha_Cd-2=\frac12$.
\end{enumerate}

\emph{Relation to previous work.} The syzygy $\rho$ is classical for unions of members of a pencil. It is the key ingredient in Vall\`es' characterization of free divisors in a pencil with smooth base locus \cite{Valles}, it is used by Dimca for curve arrangements of pencil type \cite[Section~4]{Dimca-pencils}, and it is the canonical derivation $\delta_{g,h}$ of a pencil in \cite[Lemma~3.2]{DGIMSV}. For a union of $m\ge4$ members of a pencil of curves of degree $k$, with zero-dimensional base locus and reduced members, Dimca proved that $\mdr(f)=2k-2$ \cite[Theorem~1.8 and Corollary~1.9]{Dimca-pencils}; this is the case $a=b$ of Theorem~\ref{thm:A}. The new points here are the weighted case $a\neq b$, which contains the curves of torus type, the minimal resolution in Theorem~\ref{thm:B}, and the application to Zariski pairs. In \cite{Dimca-pencils} the name \emph{pencil type} is used for unions of at least three members of a pencil; Definition~\ref{def:pencil} extends it (Remark~\ref{rk:factor}).

That Zariski pairs and Ziegler pairs are closely related was observed recently by Bannai and Tokunaga \cite{Bannai-Tokunaga}, who checked that Zariski's sextic pair has non-isomorphic Milnor algebras. Theorem~\ref{thm:C} explains the mechanism: the torus structure produces a syzygy that reaches the lowest degree allowed by \eqref{eq:DSS}, and the non-torus member stays above the bound. We propose the following converse.

\begin{conj}\label{conj:main}
If $C:f=0$ is DSS-extremal, then $C$ is of pencil type: there are forms $g,h$ with $\deg g+\deg h=\mdr(f)+2$ and a binary form $\Phi$ with $f=\Phi(g,h)$, and $\nabla g\times\nabla h$ is a syzygy of minimal degree.
\end{conj}

The conjecture holds for nodal curves (Proposition~\ref{prop:nodal}); for them we also describe $\AR(f)_{\mdr(f)}$ completely (Proposition~\ref{prop:nodal-basis}). We support it with the computations in Section~\ref{sec:comp}, done with \textsc{Singular} \cite{Singular}.

\section{Preliminaries}

\subsection{Syzygies and the du Plessis--Wall bounds}
For any $f$ the three \emph{Koszul syzygies}
\[
K_1=(0,f_z,-f_y),\quad K_2=(-f_z,0,f_x),\quad K_3=(f_y,-f_x,0),
\]
of degree $d-1$, belong to $\AR(f)$. In vector notation, $K(v)=v\times\nabla f$ for $v\in S^3$, so $K_i=e_i\times \nabla f$. Hence $\mdr(f)\le d-1$, with equality for smooth curves. If $r=\mdr(f)$, du Plessis and Wall \cite{dPW} proved
\begin{equation}\label{eq:dPW}
\tau(C)\le (d-1)^2-r(d-r-1) \ \ \text{for } r<\tfrac d2,\qquad
\tau(C)\le (d-1)^2-r(d-r-1)-\binom{2r-d+2}{2}\ \ \text{for } r\ge \tfrac d2,
\end{equation}
where $\tau(C)$ is the total Tjurina number. Equality in the first case characterizes free curves.

\subsection{Weighted homogeneous singularities}
A germ $(C,p)$ is weighted homogeneous if, in suitable local analytic coordinates, it is defined by a polynomial $\Phi(u,v)$ with $\Phi(t^{w_1'}u,t^{w_2'}v)=t^{N}\Phi(u,v)$; its weights are $w_i=w_i'/N$. For a weighted homogeneous isolated singularity with $0<w_1,w_2\le \frac12$, the Arnold exponent is $\alpha_p=w_1+w_2$ \cite[(2.1)]{DSe}. The Milnor number is $\mu=(1/w_1-1)(1/w_2-1)$ by Milnor and Orlik \cite{MO}, and $\tau=\mu$ by Saito \cite{Saito}. For instance $\alpha(A_k)=\frac12+\frac1{k+1}$ and $\alpha(D_4)=\frac23$.

\subsection{Exactness of complexes}
We use the criterion of Buchsbaum and Eisenbud \cite{BE}. A complex of finite free modules over a Noetherian ring
\[
0\to F_s\xrightarrow{\varphi_s}F_{s-1}\to\cdots\to F_1\xrightarrow{\varphi_1}F_0
\]
is exact if and only if, for every $k$, $\operatorname{rank}\varphi_{k+1}+\operatorname{rank}\varphi_k=\operatorname{rank}F_k$ and the ideal $I_{r_k}(\varphi_k)$ of $r_k\times r_k$ minors, $r_k=\operatorname{rank}\varphi_k$, has grade at least $k$. Over the Cohen--Macaulay ring $S$, the grade of a homogeneous ideal equals its codimension $3-\dim V(I)$, where $V(I)\subset\C^3$ is its zero set.

\section{Pencil-type curves are DSS-extremal}

\begin{defn}\label{def:pencil}
Let $a,b,d$ be positive integers and let $\Phi(u,v)\in\C[u,v]$ be a nonzero polynomial with $\Phi(t^au,t^bv)=t^d\Phi(u,v)$, which is reduced and belongs to $(u,v)^2$. Let $g\in S_a$, $h\in S_b$ be algebraically independent forms. If $f=\Phi(g,h)$ is reduced, we say that $C:f=0$ is a curve of \emph{pencil type} $(a,b)$. Its \emph{base locus} is $B=\{g=h=0\}\subset\PP^2$, and we set
\[
\rho=\nabla g\times\nabla h,\qquad f_g=\Phi_u(g,h),\qquad f_h=\Phi_v(g,h).
\]
\end{defn}

\begin{rk}\label{rk:factor}
Let $k=\gcd(a,b)$, $a=kA$, $b=kB$. Every monomial $u^iv^j$ of $\Phi$ satisfies $ai+bj=d$, so two exponents of $u$ differ by a multiple of $B$ and two exponents of $v$ by a multiple of $A$. Factoring out the lowest powers of $u$ and $v$, one is left with a homogeneous polynomial in $u^B$ and $v^A$, which splits into linear factors. Since $\Phi$ is reduced,
\[
\Phi=c\,u^{\varepsilon}v^{\eta}\prod_{i=1}^{m}\bigl(u^B-s_iv^A\bigr),\qquad \varepsilon,\eta\in\{0,1\},\ c,s_i\in\C^*,\ s_i\ \text{distinct}.
\]
Hence $C$ is the union of $m$ members $g^B-s_ih^A=0$ of the pencil generated by $g^B$ and $h^A$, together with $g=0$ if $\varepsilon=1$ and $h=0$ if $\eta=1$. Such pencils were used in \cite[Section~3]{DGIMSV} to construct free curves; here we look at the opposite situation, where $C$ is not free (Theorem~\ref{thm:B}). For $a=b$ we get $A=B=1$, and $C$ is a union of members of the pencil spanned by $g$ and $h$, as in \cite{Dimca-pencils}.
\end{rk}

The main examples are the following.
\begin{itemize}
\item Curves of $(p,q)$-torus type, $f=g^p+h^q$, with $\Phi=u^p+v^q$, $a=d/p$, $b=d/q$.
\item Unions of $m\ge2$ distinct members of a pencil, $f=\prod_{i=1}^m(\lambda_ig+\mu_ih)$, with $a=b$ and $d=ma$.
\item Reducible curves $f=gh$, with $\Phi=uv$.
\end{itemize}

\begin{lem}\label{lem:rho}
If $C:f=\Phi(g,h)=0$ is of pencil type $(a,b)$, then $\nabla f=f_g\nabla g+f_h\nabla h$, and $\rho$ is a nonzero element of $\AR(f)_{a+b-2}$. In particular $\mdr(f)\le a+b-2$.
\end{lem}

This is well known, see \cite{Valles}, \cite[Section~4.1]{Dimca-pencils} and \cite[Lemma~3.2]{DGIMSV}; we include the short proof for completeness.

\begin{proof}
The formula for $\nabla f$ is the chain rule. Since $\rho\cdot\nabla g=\det(\nabla g,\nabla g,\nabla h)=0$ and similarly $\rho\cdot\nabla h=0$, we get $\rho\cdot\nabla f=0$. The components of $\rho$ are $2\times 2$ minors of the Jacobian matrix of $(g,h)$, of degree $(a-1)+(b-1)$. If $\rho=0$, the differentials $dg$ and $dh$ are proportional at every point, so $g$ and $h$ are algebraically dependent by the Jacobian criterion in characteristic zero, a contradiction.
\end{proof}

\begin{lem}\label{lem:local}
Let $C$ be of pencil type $(a,b)$ and let $p\in B$ with $\rho(p)\ne0$. Then $(C,p)$ is analytically isomorphic to the germ $(\{\Phi=0\},0)$. This germ is a weighted homogeneous isolated singularity with Arnold exponent $\alpha_p=(a+b)/d$ and Milnor number $\mu_p=(d/a-1)(d/b-1)$.
\end{lem}

\begin{proof}
By the Euler formula, $p\cdot\nabla g(p)=a\,g(p)=0$ and $p\cdot\nabla h(p)=0$. In an affine chart containing $p$, say $z=1$, suppose the affine differentials of $g$ and $h$ at $p$ were proportional. Then some combination $\lambda\nabla g(p)+\mu\nabla h(p)\ne 0$ would be of the form $(0,0,c)$. Being orthogonal to $p$, it would satisfy $c=0$, which is impossible. Hence $u=g$, $v=h$ are local coordinates at $p$, and in them $f=\Phi(u,v)$.

$\Phi$ is weighted homogeneous with weights $(a/d,b/d)$. Being reduced and in $(u,v)^2$, it has an isolated singularity at the origin. If $2a\le d$ and $2b\le d$, the formula $\alpha_p=w_1+w_2$ applies and gives $(a+b)/d$; the Milnor--Orlik formula gives $\mu_p$. Suppose now $2a>d$. Then every monomial $u^iv^j$ of $\Phi$ has $i\le 1$, so $\Phi=c_1uv^j+c_2v^k$. Reducedness and $\Phi\in(u,v)^2$ force $j=1$ and $c_1\ne0$, hence $\Phi=v(c_1u+c_2v^{k-1})$. This germ consists of two smooth transversal branches, i.e.\ a node. Then $a+b=d$, $\alpha_p=1=(a+b)/d$ and $\mu_p=1=(d/a-1)(d/b-1)$. The case $2b>d$ is symmetric.
\end{proof}

The base points need not be transversal: the curves $g=0$ and $h=0$ may be tangent or singular at $p\in B$, and then the germ $(C,p)$ is no longer isomorphic to $\{\Phi=0\}$. Nevertheless, its Arnold exponent is never larger than that of $\Phi$.

\begin{prop}\label{prop:pencil-lct}
Let $C:f=\Phi(g,h)=0$ be of pencil type $(a,b)$, with arbitrary singularities. Then $\alpha_p\le (a+b)/d$ for every $p\in B$. In particular $\alpha_C\le(a+b)/d$, and the syzygy $\rho$ satisfies
\[
\deg\rho=a+b-2\ \ge\ \alpha_C\,d-2 .
\]
\end{prop}

\begin{proof}
If $g$ and $h$ had a common factor $\ell$, then $\ell^2$ would divide $f$, since $\Phi\in(u,v)^2$; this contradicts the reducedness of $f$. Hence $B$ is finite, and it is nonempty by B\'ezout. Let $p\in B$ and work in an affine chart containing $p$. Then $G=(g,h)\colon X\to Y$, where $X=(\C^2,p)$ and $Y=(\C^2,0)$, is a finite map germ, and $(C,p)$ is the pull-back $G^*D$ of the germ $D=\{\Phi=0\}$, since $f=\Phi\circ G$. Write $K_X=G^*K_Y+R$, where $R\ge0$ is the ramification divisor. Let $c=\alpha_p=\operatorname{lct}_p(C)$. Then the pair $(X,cG^*D)$ is log canonical near $p$, hence so is $(X,cG^*D-R)$. Since $K_X+cG^*D-R=G^*(K_Y+cD)$, \cite[Proposition~5.20]{KM} shows that $(Y,cD)$ is log canonical near $0$, i.e.\ $c\le\operatorname{lct}_0(\Phi)$. By Lemma~\ref{lem:local} and its proof, the germ $\{\Phi=0\}$ has Arnold exponent $(a+b)/d$, so $\alpha_p\le(a+b)/d$. The last claim follows from Lemma~\ref{lem:rho}.
\end{proof}

\begin{thm}\label{thm:A}
Let $C:f=\Phi(g,h)=0$ be of pencil type $(a,b)$ with only weighted homogeneous singularities. Then $\alpha_C\le (a+b)/d$ and $\mdr(f)\le a+b-2$. If moreover $\alpha_C=(a+b)/d$, then
\[
\mdr(f)=a+b-2=\alpha_C\,d-2,
\]
so $C$ is DSS-extremal and $\rho$ is a Jacobian syzygy of minimal degree. The equality $\alpha_C=(a+b)/d$ holds, for instance, when $\alpha_C$ is attained at a point $p\in B$ with $\rho(p)\neq0$.
\end{thm}

\begin{proof}
Proposition~\ref{prop:pencil-lct} gives $\alpha_C\le(a+b)/d$, and Lemma~\ref{lem:rho} gives $\mdr(f)\le a+b-2$. If $\alpha_C=(a+b)/d$, the bound \eqref{eq:DSS} gives $\mdr(f)\ge a+b-2$. If $\alpha_C=\alpha_p$ for some $p\in B$ with $\rho(p)\ne0$, then $\alpha_C=(a+b)/d$ by Lemma~\ref{lem:local}.
\end{proof}

\begin{rk}\label{rk:m3}
For a union of three members of a pencil of curves of degree $k$, Dimca found cases with $\mdr(f)<2k-2$, some of them free \cite[Theorem~1.8]{Dimca-pencils}. This does not contradict Theorem~\ref{thm:A}: here $(a+b)/d=2/3$, and by \eqref{eq:DSS} such a curve must have a singularity that is not weighted homogeneous, or a singularity with Arnold exponent smaller than $2/3$, for instance at a base point where $\rho$ vanishes.
\end{rk}

\begin{cor}\label{cor:torus}
Let $C:g^p+h^q=0$ be a curve of $(p,q)$-torus type, $p,q\ge2$, of degree $d$, all of whose singular points are base points with $\rho\ne0$. Then $C$ has $(d/p)(d/q)$ singular points of type $u^p+v^q=0$, $\alpha_C=\frac1p+\frac1q$, and
\[
\mdr(f)=\frac dp+\frac dq-2.
\]
\end{cor}

For $(p,q)=(2,3)$ and $d=6$ this is Zariski's sextic with six cusps on the conic $h=0$, with $\mdr(f)=3$. For $d=6m$ one obtains $\mdr(f)=5m-2$.

\begin{prop}\label{prop:nodal}
Conjecture~\ref{conj:main} holds for nodal curves.
\end{prop}

\begin{proof}
For a nodal curve $\alpha_C=1$, and Dimca and Sernesi showed that $\mdr(f)=d-2$ holds exactly when $C$ is reducible \cite[Example 2.6]{DSe}. If $C=C_1\cup C_2$ with $C_1:g=0$ and $C_2:h=0$, then $f=gh$ is of pencil type with $\Phi=uv$, and by Lemma~\ref{lem:rho} the syzygy $\nabla g\times\nabla h$ has degree $\deg g+\deg h-2=d-2=\mdr(f)$.
\end{proof}

\begin{prop}\label{prop:nodal-basis}
Let $C:f=0$ be a nodal curve with irreducible components $C_j:f_j=0$, $j=1,\dots,r$, $r\ge2$. Each splitting $f=f_j\cdot(f/f_j)$ is a pencil structure of type $\Phi=uv$, with syzygy $\rho_j=\nabla f_j\times\nabla(f/f_j)\in\AR(f)_{d-2}$. Then $\rho_1+\dots+\rho_r=0$, and $\rho_1,\dots,\rho_{r-1}$ form a basis of $\AR(f)_{d-2}$. In particular $\AR(f)_{\mdr(f)}$ is spanned by the vectors $\nabla g\times\nabla h$ over all pencil structures of $C$.
\end{prop}

\begin{proof}
Set $\omega_{jk}=\frac{f}{f_jf_k}\,\nabla f_j\times\nabla f_k$ for $j\neq k$, so that $\omega_{kj}=-\omega_{jk}$. Since $\nabla(f/f_j)=\sum_{k\neq j}\frac{f}{f_jf_k}\nabla f_k$, we get $\rho_j=\sum_{k\neq j}\omega_{jk}$, and the antisymmetry gives $\sum_j\rho_j=0$.

Suppose $\sum_{j=1}^{r}c_j\rho_j=0$ with $c_j\in\C$. Let $l\neq k$ and let $p\in C_l\cap C_k$, which is nonempty by B\'ezout. Since $C$ is nodal, $p$ lies on no other component, and $C_l$, $C_k$ are smooth at $p$ with distinct tangents. Every $\omega_{ab}$ with $\{a,b\}\neq\{l,k\}$ contains the factor $f_l$ or $f_k$, so it vanishes at $p$. Hence
\[
0=\sum_jc_j\rho_j(p)=(c_l-c_k)\,\omega_{lk}(p),
\]
and $\omega_{lk}(p)\neq0$ because $\nabla f_l(p)$, $\nabla f_k(p)$ are independent and $f/(f_lf_k)$ does not vanish at $p$. So $c_l=c_k$ for all $l,k$: the only relations among $\rho_1,\dots,\rho_r$ are multiples of $\sum_j\rho_j=0$, and $\rho_1,\dots,\rho_{r-1}$ are linearly independent. By \cite[Theorem~4.1]{DS-Edinburgh}, $\dim\AR(f)_{d-2}=r-1$, which ends the proof.
\end{proof}

This answers Question~(2) of Section~\ref{sec:final} for nodal curves. We also checked it in \textsc{Singular} for a random union of two conics and a line and of a cubic, a conic and two lines (Section~\ref{sec:comp}).

\section{The minimal resolution of a transversal pencil-type curve}

We now impose the following transversality hypothesis:
\[
\textbf{(H)}\qquad \rho(p)\ne 0\ \text{ for every } p\in\Sing(C).
\]

\begin{lem}\label{lem:H}
Let $C:f=\Phi(g,h)=0$ be of pencil type $(a,b)$ satisfying \textbf{(H)}. Then:
\begin{enumerate}[(i)]
\item $\Sing(C)=B$, and the curves $g=0$ and $h=0$ meet transversally in $ab$ points;
\item every singular point of $C$ is analytically of type $\Phi$, $\alpha_C=(a+b)/d$, and $C$ is DSS-extremal;
\item $\tau(C)=\mu(C)=(d-a)(d-b)$.
\end{enumerate}
\end{lem}

\begin{proof}
Let $p\in\Sing(C)$. Then $0=\nabla f(p)=f_g(p)\nabla g(p)+f_h(p)\nabla h(p)$ with $\nabla g(p),\nabla h(p)$ independent, hence $\Phi_u(g(p),h(p))=\Phi_v(g(p),h(p))=0$. The common zero set of $\Phi_u$ and $\Phi_v$ in $\C^2$ is invariant under the action $t\cdot(u,v)=(t^au,t^bv)$, and every orbit has the origin in its closure. Since the origin is an isolated singular point of $\Phi$, the only common zero is $(0,0)$, and so $p\in B$. Conversely, every point of $B$ is singular. At every $p\in B$ we have $\rho(p)\neq 0$, so $g=0$ and $h=0$ are smooth at $p$ with distinct tangents; the intersection is transversal and $|B|=ab$ by B\'ezout. This proves (i). Part (ii) follows from Lemma~\ref{lem:local} and Theorem~\ref{thm:A}. For (iii), $\tau(C)=\sum_p\tau_p=ab\,\mu_p=ab(d/a-1)(d/b-1)=(d-a)(d-b)$.
\end{proof}

\begin{thm}\label{thm:B}
Let $C:f=\Phi(g,h)=0$ be of pencil type $(a,b)$ with $a,b\ge2$, satisfying \textbf{(H)}. Then $\AR(f)$ is minimally generated by $\rho,K_1,K_2,K_3$. The only minimal relations among them are
\begin{equation}\label{eq:rel}
R_1:\ f_g\,\rho+\sum_{i}h_{x_i}K_i=0,\qquad R_2:\ -f_h\,\rho+\sum_i g_{x_i}K_i=0 .
\end{equation}
Consequently $C$ is a $4$-syzygy curve with exponents $(a+b-2,\,d-1,\,d-1,\,d-1)$, and $M(f)$ has the minimal graded free resolution
\[
0\to S(-2d+3-a)\oplus S(-2d+3-b)\xrightarrow{\varphi_3}
S(-d-a-b+3)\oplus S(-2d+2)^3\xrightarrow{\varphi_2}
S(1-d)^3\xrightarrow{\varphi_1}S\to M(f)\to0 .
\]
\end{thm}

\begin{proof}
Write $\sum_ih_{x_i}K_i=\nabla h\times\nabla f$. By Lemma~\ref{lem:rho},
\[
\nabla h\times\nabla f=f_g\,\nabla h\times\nabla g=-f_g\,\rho,\qquad
\nabla g\times\nabla f=f_h\,\nabla g\times\nabla h=f_h\,\rho,
\]
which gives \eqref{eq:rel}. Let $\varphi_1=(f_x,f_y,f_z)$. Let $\varphi_2$ be the $3\times 4$ matrix with columns $\rho,K_1,K_2,K_3$. Let $\varphi_3$ be the $4\times2$ matrix with columns $(f_g,h_x,h_y,h_z)^T$ and $(-f_h,g_x,g_y,g_z)^T$. Then $\varphi_1\varphi_2=0$ by Lemma~\ref{lem:rho}, and $\varphi_2\varphi_3=0$ by \eqref{eq:rel}. With the indicated shifts the maps are homogeneous of degree zero.

\emph{Ranks.} We have $\operatorname{rank}\varphi_1=1$. At a point $q$ with $\nabla f(q)\ne0$, the columns of $\varphi_2$ lie in the plane $\nabla f(q)^{\perp}$, and the skew-symmetric matrix $(K_1,K_2,K_3)$ has rank $2$; hence $\operatorname{rank}\varphi_2=2$. The $2\times2$ minors of $\varphi_3$ built from rows $i+1,j+1$ are, up to sign, the components of $\rho\ne0$; hence $\operatorname{rank}\varphi_3=2$. The rank conditions $2+2=4$, $2+1=3$ hold.

\emph{Grades.} $I_1(\varphi_1)=J_f\ne0$ has grade $\ge1$. If $q\in\C^3$ satisfies $\nabla f(q)\ne0$, some $2\times2$ minor of $(K_1,K_2,K_3)$ is nonzero at $q$, so $V(I_2(\varphi_2))\subset V(J_f)$. The latter is the affine cone over $\Sing(C)$, of dimension $\le1$, so the grade of $I_2(\varphi_2)$ is at least $2$. Finally, the minor of $\varphi_3$ on rows $1$ and $i+1$ is $f_g g_{x_i}+f_h h_{x_i}=f_{x_i}$, and the minors on rows $i+1,j+1$ are $\pm\rho_k$. Hence
\[
I_2(\varphi_3)\supseteq J_f+(\rho_1,\rho_2,\rho_3).
\]
If $0\ne q\in V(J_f)\cap V(\rho)$, then $f(q)=\frac1d\,q\cdot\nabla f(q)=0$, so $[q]\in\Sing(C)$ and $\rho([q])=0$, contradicting \textbf{(H)}. Thus $V(I_2(\varphi_3))=\{0\}$ and its grade is $3$.

By the Buchsbaum--Eisenbud criterion the complex is exact, so it is a free resolution of $M(f)$. In particular $\AR(f)=\ker\varphi_1=\operatorname{im}\varphi_2$ is generated by $\rho,K_1,K_2,K_3$, and $\ker\varphi_2$ is generated by $R_1,R_2$. The resolution is minimal because no matrix entry is a nonzero constant. The entries of $\varphi_2$ have degrees $a+b-2\ge2$ and $d-1$. The entries of $\varphi_3$ have degrees $d-a$, $d-b$, $a-1$, $b-1$, all positive since $a,b\ge2$ and $d>a,b$. The last inequalities hold because $\Phi\in(u,v)^2$ is reduced, so both $u$ and $v$ occur in $\Phi$ and no monomial of $\Phi$ is linear. The twists follow from the degrees: $\deg\rho=a+b-2$, $\deg K_i=d-1$, and $R_1,R_2$ have degrees $b-1+d-1$ and $a-1+d-1$ as elements of $S^4$ graded by the generators.
\end{proof}

\begin{rk}\label{rk:betti}
Theorem~\ref{thm:B} shows that all transversal pencil-type curves with the same $(d,a,b)$ have the same graded Betti numbers, whatever $\Phi$ is. For instance, for $(p,q)$-torus curves with $d=6$, $(p,q)=(2,3)$, the Betti table has $\AR(f)$ generated in degrees $3,5,5,5$ with two second syzygies in degrees $6$ and $7$. This is exactly what \textsc{Singular} returns for Zariski's sextic and for random torus sextics.
\end{rk}

\begin{rk}\label{rk:uv}
Let $\Phi=uv$, so $f=gh$ and $d=a+b$. A singular point of $g=0$ lies on $C$, is singular for $C$, and there $\nabla g=0$, hence $\rho=0$. So \textbf{(H)} forces $g=0$ and $h=0$ to be smooth, hence irreducible, and $C$ has exactly two components meeting transversally in $ab$ nodes. In this case Theorem~\ref{thm:B} agrees with \cite[Theorem~4.1]{DS-Edinburgh}, which gives $r-1=1$ syzygy of degree $d-2$, and with the minimal resolution of the Jacobian ideal of a product of general forms in \cite[Example~3.7]{BRST} for two factors. (There the twist $R(-2D+d_j-1)$ should read $R(-2D-d_j+1)$; with this correction the two resolutions coincide.) If a nodal curve has $r\ge3$ components, then every splitting $f=gh$ has a reducible factor and \textbf{(H)} fails. For instance $f=(x^3+y^3+z^3)\,xy$ is of pencil type $(3,2)$ with $\Phi=uv$, \textbf{(H)} fails at $[0:0:1]$, and $\dim\AR(f)_3=r-1=2$ by \cite[Theorem~4.1]{DS-Edinburgh}, so Theorem~\ref{thm:B} cannot hold for it. The hypothesis $a,b\ge2$ excludes the case where $g$ or $h$ is a line, in which the resolution of Theorem~\ref{thm:B} is not minimal.
\end{rk}

\begin{ex}\label{ex:H-fails}
Hypothesis \textbf{(H)} cannot be dropped.
\begin{enumerate}[(1)]
\item The arrangement $\mathcal A:(x^2-y^2)(y^2-z^2)(z^2-x^2)=0$ is of pencil type $(2,2)$, since $x^2-y^2$, $y^2-z^2$, $z^2-x^2$ belong to one pencil. It is DSS-extremal, with $\mdr=2$. Besides the four triple points in $B$ it has three nodes, at which a member of the pencil is singular and $\rho=0$. The arrangement is the union of the three singular members of the pencil spanned by $x^2-y^2$ and $y^2-z^2$, whose base locus consists of four points. By Vall\`es' theorem \cite{Valles} (see also \cite[Theorem~1.5]{DGIMSV} and \cite[Theorem~1.14]{Dimca-pencils}) it is free with exponents $(2,3)$, and $\tau=19$, so Theorem~\ref{thm:B} does not apply.
\item The sextic $x^6+y^6+z^6-2(x^3y^3+x^3z^3+y^3z^3)=0$ with nine cusps is of torus type. Degtyarev showed it has twelve torus structures \cite{Degtyarev}. For each torus structure at most six of the nine cusps lie in its base locus, so \textbf{(H)} fails. The curve is DSS-extremal ($\mdr=3$) but nearly free with exponents $(3,3)$, so $\dim\AR(f)_3=3$: the twelve vectors $\nabla g\times\nabla h$ span a $3$-dimensional space.
\end{enumerate}
\end{ex}

Several classical line arrangements are of pencil type. Theorem~\ref{thm:A} shows that some of them are DSS-extremal, while for others the condition $\alpha_C=(a+b)/d$ fails. An ordinary point of multiplicity $m\ge2$ has Arnold exponent $2/m$; in particular $\alpha(D_4)=\frac23$.

\begin{ex}\label{ex:arr}
\begin{enumerate}[(1)]
\item \emph{The Ceva arrangement} $\mathcal A(3,3,3):\ f=(x^3-y^3)(y^3-z^3)(z^3-x^3)=0$, of degree $d=9$. It is the union of three members of the pencil of cubics spanned by $g=x^3-y^3$ and $h=y^3-z^3$, so it is of pencil type $(3,3)$ with $(a+b)/d=\frac23$. Its $12$ singular points are ordinary triple points, so $\alpha_C=\frac23$, and Theorem~\ref{thm:A} gives $\mdr(f)=4=\alpha_Cd-2$. As for $\mathcal A(2,2,3)$, hypothesis \textbf{(H)} fails at the three points $[1:0:0]$, $[0:1:0]$, $[0:0:1]$, where a member of the pencil is singular. The arrangement is free with exponents $(4,4)$ and $\tau=48$, so $\dim\AR(f)_4=2$. The nine lines can also be grouped into three triangles belonging to one pencil of cubics, and \textsc{Singular} shows that this happens in exactly three ways. Together with the pencil above, this gives four pencil structures of type $(3,3)$, and their four vectors $\nabla g\times\nabla h$ span $\AR(f)_4$. This supports Question~(2) of Section~\ref{sec:final} beyond nodal curves.
\item \emph{The Hesse arrangement} of $12$ lines, formed by the four triangles of the Hesse pencil spanned by $g=x^3+y^3+z^3$ and $h=xyz$. Since $\prod_{k=0}^{2}(g-3\omega^kh)=g^3-27h^3$, where $\omega$ is a primitive cube root of unity, it is
\[
f=xyz\bigl((x^3+y^3+z^3)^3-27x^3y^3z^3\bigr)=\Phi(g,h),\qquad \Phi=v(u^3-27v^3),
\]
of pencil type $(3,3)$ with $d=12$ and $(a+b)/d=\frac12$. Its singular points are the nine base points of the pencil, which are ordinary quadruple points, and the twelve vertices of the triangles, which are nodes and at which $\rho=0$. Hence $\alpha_C=\frac12$ is attained on $B$, and Theorem~\ref{thm:A} gives $\mdr(f)=4=\alpha_Cd-2$: the Hesse arrangement is DSS-extremal. It is free with exponents $(4,7)$ and $\tau=93$, and $\AR(f)_4$ is spanned by $\rho$.
\item \emph{The Fermat arrangements} $\mathcal A(m,m,3):\ (x^m-y^m)(y^m-z^m)(z^m-x^m)=0$, $m\ge2$, of degree $d=3m$, are of pencil type $(m,m)$ with $(a+b)/d=\frac23$. For $m=2,3$ they are DSS-extremal, see Example~\ref{ex:H-fails}(1) and part~(1). For $m\ge4$ the three coordinate points are ordinary $m$-fold points, so $\alpha_C=\frac2m<\frac23$, and Theorem~\ref{thm:A} does not apply. Here $\alpha_Cd-2=4$ for every $m$, while $\mdr(f)=m+1$, since $\mathcal A(m,m,3)$ is the reflection arrangement of the group $G(m,m,3)$, which is free with exponents $(m+1,2m-2)$ \cite{OT}. So these curves are of pencil type but not DSS-extremal, and since $m+1<2m-2$, the syzygy $\rho$ is not of minimal degree; compare Remark~\ref{rk:m3}. \textsc{Singular} confirms $\mdr(f)=5,6$ and $\tau=91,148$ for $m=4,5$.
\end{enumerate}
\end{ex}

\section{Sextics with six cusps and Ziegler pairs}

Zariski showed that irreducible sextics with six ordinary cusps form at least two families: those of torus type $g^2+h^3$ ($g\in S_3$, $h\in S_2$), whose cusps lie on the conic $h=0$, and those whose cusps are not on a conic \cite{Zariski}. The two families give the first Zariski pair \cite{ACT}. An explicit non-torus sextic was given by Oka \cite{Oka} (affine coordinates):
\begin{equation}\label{eq:Oka}
x^2(x-1)^2(x^2+2x)-x^2(x-1)^2(y^2-1)+\tfrac13(x^2-1)(y^2-1)^2-\tfrac1{27}(y^2-1)^3=0 .
\end{equation}

\begin{lem}\label{lem:inner}
Let $C:g^2+h^3=0$ be a torus sextic whose only singularities are six ordinary cusps. Then $C$ satisfies \textbf{(H)}.
\end{lem}

\begin{proof}
Every point of $B$ is singular. At $p\in B$ where the conic $h=0$ and the cubic $g=0$ are smooth and meet transversally, $\rho(p)\ne0$, and the germ is $u^2+v^3$, an $A_2$ cusp. In all other cases the germ at $p$ is not an $A_2$ cusp:
\begin{itemize}
\item if the curves are tangent at $p$, the local equation is $w^2+(w+cu^2)^3$ in suitable coordinates, an $A_{\ge5}$ singularity;
\item if the cubic is singular at $p$, the germ $g^2+h^3$ has multiplicity $\ge3$;
\item if the conic is singular at $p$ and the cubic is smooth there, the germ has the form $u^2+(\text{terms of order}\ge 6)$, again not $A_2$.
\end{itemize}
Since all singularities are $A_2$, the intersection is transversal at every point of $B$, so $|B|=6$. These six cusps exhaust $\Sing(C)$, hence \textbf{(H)} holds.
\end{proof}

\begin{rk}
The three degenerate cases in the proof of Lemma~\ref{lem:inner} were also checked in \textsc{Singular}. For a cubic tangent to the conic at a base point, for a singular conic, and for a cubic singular at a base point, one obtains respectively $\tau(C)=13,13,14$ and five singular points in each case. This is consistent with the configurations $4A_2+A_5$, $4A_2+A_5$ and $4A_2+E_6$, never with $6A_2$.
\end{rk}

\begin{thm}\label{thm:C}
Let $C:f=0$ be an irreducible sextic whose singularities are exactly six ordinary cusps. Then:
\begin{enumerate}[(i)]
\item $3\le\mdr(f)\le4$;
\item if $C$ is of torus type, then $\mdr(f)=3$, $C$ is DSS-extremal, and $\AR(f)$ is minimally generated in degrees $3,5,5,5$ as in Theorem~\ref{thm:B};
\item Oka's sextic \eqref{eq:Oka} has $\mdr(f)=4$, and $\AR(f)$ is minimally generated in degrees $4,4,5,5,5$;
\item if $\Sigma$ is the irreducible component of the variety of six-cuspidal sextics containing Oka's curve, then $\mdr(f)=4$ for the general member of $\Sigma$.
\end{enumerate}
In particular, Zariski's torus sextic and Oka's sextic form a Zariski pair which is also a Ziegler pair, distinguished by $\mdr$.
\end{thm}

\begin{proof}
(i) Since $\alpha(A_2)=5/6$, \eqref{eq:DSS} gives $\mdr(f)\ge3$. Here $\tau(C)=12$. If $r=\mdr(f)=5$, the second bound in \eqref{eq:dPW} gives $\tau(C)\le 25-0-\binom62=10<12$, a contradiction. Since $\mdr(f)\le d-1=5$, we get $\mdr(f)\le4$.

(ii) By Lemma~\ref{lem:inner}, $C$ satisfies \textbf{(H)} with $(a,b)=(3,2)$. Lemma~\ref{lem:H} and Theorem~\ref{thm:B} give $\mdr(f)=3$ and the stated generators.

(iii) This is a direct computation in \textsc{Singular} (Section~\ref{sec:comp}). The homogenized equation $27\cdot$\eqref{eq:Oka} is irreducible with $\tau=12$, and the minimal resolution of $M(f)$ has $\AR(f)$ generated in degrees $4,4,5,5,5$, with three second syzygies.

(iv) The map $f\mapsto\dim\AR(f)_3=\dim\ker\bigl(S_3^3\to S_8,\ (a,b,c)\mapsto af_x+bf_y+cf_z\bigr)$ is upper semicontinuous on the space of sextics. It vanishes at Oka's curve, hence on a nonempty Zariski open subset of $\Sigma$. On that subset $\mdr(f)\ge4$, and so $\mdr(f)=4$ by (i).
\end{proof}

\begin{rk}
If the variety of non-torus six-cuspidal sextics is irreducible, Theorem~\ref{thm:C}(iv) holds for the general non-torus sextic. Oka's curve has the semi-torus form $f_2^3+g_2^2h_2$ \cite{Oka}, which is not of pencil type, in agreement with Conjecture~\ref{conj:main}.
\end{rk}

\begin{conj}\label{conj:sextic}
Every irreducible six-cuspidal sextic of non-torus type has $\mdr(f)=4$. Equivalently, an irreducible six-cuspidal sextic is DSS-extremal if and only if it is of torus type.
\end{conj}

A natural approach to Conjecture~\ref{conj:sextic} is the relation between Jacobian syzygies and defects of linear systems through the singular points \cite{Dimca-defects}. One expects $\AR(f)_3\neq0$ exactly when the six cusps impose dependent conditions on conics, i.e.\ when they lie on a conic, which by Zariski's classical result characterizes the torus type. For cuspidal curves of degree $6k$, Kloosterman related the defect of the cusps with respect to curves of degree $5k-3$ to the syzygies of the ideal of the cusps and to a Mordell--Weil rank \cite{Kloosterman}; for $k=1$ these are exactly the conics through the six cusps.

\section{Computations}\label{sec:comp}

All computations were done in \textsc{Singular} \cite{Singular}. We computed $\mdr(f)$ as the minimal degree of the generators of $\Syz(J_f)$, $\tau(C)$ as the degree of $S/J_f$, and the minimal resolution with \texttt{mres}. Random forms have integer coefficients in $[-9,9]$. For the random torus sextic, the two nonics with nine $D_4$ points and the union of a cubic and a conic we also checked, by Gr\"obner basis reduction, that $\AR(f)$ is contained in the submodule generated by $\rho,K_1,K_2,K_3$. For the two random nodal curves of Proposition~\ref{prop:nodal-basis}, with $r=3$ and $r=4$ components, \textsc{Singular} returns $\dim\AR(f)_{d-2}=r-1$, and the minimal syzygies of degree $d-2$ lie in the span of the $\rho_j$. A \textsc{Singular} script reproducing all these computations is available from the author upon request.

\begin{table}[ht]
\centering\small
\begin{tabular}{@{}llccccl@{}}
\toprule
Curve & Singularities & $\alpha_C$ & $\alpha_Cd-2$ & $\mdr$ & $\tau$ & Generators of $\AR(f)$\\
\midrule
Zariski $(x^2+y^2)^3+(y^3+z^3)^2$ & $6A_2$ & $5/6$ & $3$ & $3$ & $12$ & $3,5,5,5$\\
random $g^2+h^3$, $d=6$ & $6A_2$ & $5/6$ & $3$ & $3$ & $12$ & $3,5,5,5$\\
random $g^2+h^3$, $d=12$ & $24A_2$ & $5/6$ & $8$ & $8$ & $48$ & $8,11,11,11$\\
random $g^2+h^5$, $d=10$ & $10A_4$ & $7/10$ & $5$ & $5$ & $40$ & $5,9,9,9$\\
random $g^3+h^4$, $d=12$ & $12E_6$ & $7/12$ & $5$ & $5$ & $72$ & $5,11,11,11$\\
$(x^3+y^3+z^3)^3+(x^3+2y^3+3z^3)^3$ & $9D_4$ & $2/3$ & $4$ & $4$ & $36$ & $4,8,8,8$\\
$gh(g+h)$, $g,h$ random cubics & $9D_4$ & $2/3$ & $4$ & $4$ & $36$ & $4,8,8,8$\\
$gh$, cubic $\cdot$ conic & $6A_1$ & $1$ & $3$ & $3$ & $6$ & $3,4,4,4$\\
\midrule
$(x^2-y^2)(y^2-z^2)(z^2-x^2)$ & $4D_4+3A_1$ & $2/3$ & $2$ & $2$ & $19$ & $2,3$ (free)\\
nine-cuspidal sextic & $9A_2$ & $5/6$ & $3$ & $3$ & $18$ & $3,3,3$ (nearly free)\\
Ceva $\mathcal A(3,3,3)$ & $12D_4$ & $2/3$ & $4$ & $4$ & $48$ & $4,4$ (free)\\
Hesse arrangement & $9X_9+12A_1$ & $1/2$ & $4$ & $4$ & $93$ & $4,7$ (free)\\
\midrule
Oka's non-torus sextic & $6A_2$ & $5/6$ & $3$ & $4$ & $12$ & $4,4,5,5,5$\\
Fermat $\mathcal A(4,4,3)$ & $3X_9+16D_4$ & $1/2$ & $4$ & $5$ & $91$ & $5,6$ (free)\\
Fermat $\mathcal A(5,5,3)$ & $3O_5+25D_4$ & $2/5$ & $4$ & $6$ & $148$ & $6,8$ (free)\\
\bottomrule
\end{tabular}
\caption{Extremality and Jacobian syzygies. The first block satisfies \textbf{(H)}; there $\tau=(d-a)(d-b)$ and the resolution is the one of Theorem~\ref{thm:B}, including the two second syzygies in the predicted degrees. The second block is of pencil type without \textbf{(H)}. The last block is not DSS-extremal; the Fermat arrangements are of pencil type, but $\alpha_C<(a+b)/d$. Here $X_9$ denotes an ordinary quadruple point and $O_5$ an ordinary quintuple point.}
\label{tab:comp}
\end{table}

\section{Beyond weighted homogeneous singularities}\label{sec:nonWH}

The proof of \eqref{eq:DSS} in \cite{DSa} uses the weighted homogeneity of the singularities in an essential way: by \cite[Theorem~5.2]{DSa}, if some singularity of $C$ is not weighted homogeneous, then the pole order spectral sequence need not degenerate at $E_2$. For an arbitrary isolated plane curve singularity $(C,p)$ we let $\alpha_p$ be its Arnold exponent, i.e.\ the minimal spectral number, which for plane curves coincides with the log canonical threshold $\operatorname{lct}_p(C)$, and we set $\alpha_C=\min_p\alpha_p$. For weighted homogeneous germs this is the invariant used above. We do not know any curve violating \eqref{eq:DSS}, and we propose the following.

\begin{conj}\label{conj:nonWH}
Every reduced plane curve $C:f=0$ of degree $d$ satisfies $\mdr(f)\ge\alpha_C\,d-2$.
\end{conj}

A counterexample must have a singularity which is not weighted homogeneous, and such a singularity is far from being mild. Moreover, by Proposition~\ref{prop:pencil-lct}, which holds without any assumption on the singularities, a Jacobian syzygy of minimal degree of a counterexample cannot be of the form $\nabla g\times\nabla h$ for a pencil structure $f=\Phi(g,h)$.

\begin{lem}\label{lem:half}
If $(C,p)$ is not weighted homogeneous, then $\alpha_p\le\frac12$. Hence a curve violating \eqref{eq:DSS} has $\mdr(f)<\frac d2-2$.
\end{lem}

\begin{proof}
Simple singularities are weighted homogeneous, so $(C,p)$ is not simple. If $\operatorname{mult}_pC\ge4$, then $\alpha_p\le 2/\operatorname{mult}_pC\le\frac12$, using the semicontinuity of Arnold exponent \cite{Ku,V8}. Otherwise $\operatorname{mult}_pC=3$, the $3$-jet is a cube $u^3$, and by Arnold's classification there are no terms $v^4$, $uv^3$, $v^5$ as for $E_6$, $E_7$, $E_8$. In suitable coordinates every monomial of the local equation then has weighted order at least $1$ for the weights $(\frac13,\frac16)$, and the monomial valuation bound gives $\alpha_p\le\frac13+\frac16=\frac12$. The second claim follows since $\mdr(f)<\alpha_Cd-2\le\frac d2-2$.
\end{proof}

The main result of this section proves Conjecture~\ref{conj:nonWH} for curves with a Jacobian syzygy of degree at most $1$. These curves were classified by du Plessis and Wall \cite{dPW99}: a reduced curve has a syzygy of degree $1$ if and only if it is invariant under a one-parameter group of projective transformations \cite[Proposition~1.1]{dPW99}. The case of a unipotent group gives an infinite family with a singularity which is not weighted homogeneous. For $k\ge1$, $\varepsilon\in\{0,1\}$ and distinct $c_1,\dots,c_k\in\C$ consider
\begin{equation}\label{eq:Cke}
C_{k,\varepsilon}:\ z^{\varepsilon}\prod_{i=1}^{k}\bigl(y^2-2xz-c_iz^2\bigr)=0,\qquad d=2k+\varepsilon .
\end{equation}
Up to the exchange of $x$ and $z$, these are the curves with symmetry group $\mathit{Nil}_3$ in \cite[Section~2]{dPW99}: $k$ smooth conics meeting in a single point, plus possibly their common tangent line at this point.

\begin{prop}\label{prop:mdr1}
Every reduced plane curve with $\mdr(f)\le1$ satisfies $\mdr(f)\ge\alpha_Cd-2$. Up to projective equivalence, those having a singularity which is not weighted homogeneous are exactly the curves $C_{k,\varepsilon}$ with $k\ge3$ if $\varepsilon=0$ and $k\ge2$ if $\varepsilon=1$. For them
\[
\mdr(f)=1,\qquad \alpha_C=\frac{5}{4k+2\varepsilon},\qquad \alpha_C\,d-2=\frac12 .
\]
\end{prop}

\begin{proof}
If $\mdr(f)=0$, then $C$ consists of $d$ concurrent lines; its only singularity is an ordinary $d$-fold point, with $\alpha_C=2/d$, so $\alpha_Cd-2=0$.

Let $\mdr(f)=1$. A syzygy of degree $1$ is given by a matrix $A\in\mathfrak{gl}_3(\C)$ with $(Av)\cdot\nabla f(v)=0$, i.e.\ $f$ is invariant under the flow $\exp(tA)$. We recall the cases of \cite[Section~2]{dPW99}. If $A$ is diagonalizable, $f$ is semi-invariant under a torus, its singular points are fixed points of the torus, and there the local equations are weighted homogeneous; so \eqref{eq:DSS} holds by \cite{DSa,DSe}. If $A$ is conjugate to $J_2(a)\oplus(b)$, the decomposition of $f$ into weight spaces forces $f\in\C[y,z]$ in suitable coordinates, hence $\mdr(f)=0$. If $A$ is conjugate to $J_3(a)$, invariance forces $a=0$; the invariants of $y\partial_x+z\partial_y$ are generated by $z$ and $y^2-2xz$, so $f$ is as in \eqref{eq:Cke} up to a constant.

For $C_{k,\varepsilon}$, the only singular point is $p=[1:0:0]$ \cite[Proposition~2.2]{dPW99}: any two of the conics have intersection multiplicity $4$ at $p$, and the line $z=0$ is tangent to each of them at $p$. In the affine chart $x=1$ put $w=2z-y^2$. The branches of the conics are $w=-\frac{c_i}{4}y^4+\cdots$, and the line is $w+y^2=0$. Hence the Newton polygon of the local equation in the coordinates $(w,y)$ has vertices $(k,0)$ and $(0,4k)$ if $\varepsilon=0$, and $(k+1,0)$, $(k,2)$, $(0,4k+2)$ if $\varepsilon=1$. For $k\ge2$ the diagonal meets the edge $4i+j=4k+2\varepsilon$ in the point $(t,t)$ with $t=(4k+2\varepsilon)/5$. The polynomial of this edge is $y^{2\varepsilon}\prod_i\bigl(w+\frac{c_i}{4}y^4\bigr)$, which is square-free outside the coordinate axes since the $c_i$ are distinct. So the germ is Newton non-degenerate on this edge, and $\alpha_p=1/t=5/(4k+2\varepsilon)$ \cite{Varchenko}. Therefore
\[
\alpha_C\,d-2=\frac{5(2k+\varepsilon)}{2(2k+\varepsilon)}-2=\frac12<1=\mdr(f).
\]
Finally, $\mu_p=(k-1)(4k-1)$ for $\varepsilon=0$ and $\mu_p=k(4k-1)$ for $\varepsilon=1$, while $\tau(C)=\tau_p=d^2-3d+3$ \cite[Proposition~3.1]{dPW99}. So $\mu_p>\tau_p$, i.e.\ the singularity is not weighted homogeneous, exactly when $d\ge5$. For $(k,\varepsilon)=(2,0)$ the singularity is of type $A_7$, and for $k=1$ the curve is a smooth conic or a conic with a tangent line (an $A_3$ point).
\end{proof}

\begin{rk}\label{rk:mdr1-free}
The classification and the values of $\mu_p$ and $\tau(C)$ are due to du Plessis and Wall \cite{dPW99}; what is new in Proposition~\ref{prop:mdr1} is the value of $\alpha_C$ and the fact that $\alpha_Cd-2=\frac12$ for the whole family. Since $\tau(C)=(d-1)^2-(d-2)$ is the maximal value allowed by \eqref{eq:dPW} for $r=1$, the curves $C_{k,\varepsilon}$ are free with exponents $(1,d-2)$, although their singularity is not weighted homogeneous; this is also contained in the proof of \cite[Proposition~3.1]{dPW99}. Table~\ref{tab:nonWH} gives the invariants for small $k$.
\end{rk}

\begin{table}[ht]
\centering\small
\begin{tabular}{@{}cccccccc@{}}
\toprule
$k$ & $\varepsilon$ & $d$ & $\mdr$ & $\tau(C)$ & $(\mu_p,\tau_p)$ & $\alpha_C$ & $\alpha_Cd-2$\\
\midrule
2&0&4&1&7&$(7,7)$, $A_7$&$5/8$&$1/2$\\
2&1&5&1&13&$(14,13)$&$1/2$&$1/2$\\
3&0&6&1&21&$(22,21)$&$5/12$&$1/2$\\
3&1&7&1&31&$(33,31)$&$5/14$&$1/2$\\
4&0&8&1&43&$(45,43)$&$5/16$&$1/2$\\
5&0&10&1&73&$(76,73)$&$1/4$&$1/2$\\
6&1&13&1&133&$(138,133)$&$5/26$&$1/2$\\
\bottomrule
\end{tabular}
\caption{The curves $C_{k,\varepsilon}$. Since $\tau(C)=\tau_p$, the point $p=[1:0:0]$ is the only singular point, and $\mu_p>\tau_p$ means that it is not weighted homogeneous. The values of $\alpha_C$ were also checked with the \texttt{spectrum} procedure of \textsc{Singular}.}
\label{tab:nonWH}
\end{table}

\begin{rk}\label{rk:add-curve}
Adding a curve to $C_{k,\varepsilon}$ produces more examples. If $C:f=0$ has a Jacobian syzygy of degree $r$ and $D:g=0$ is a curve of degree $\beta$ with no common component with $C$, then $fg$ has a Jacobian syzygy of degree $r+\beta$ \cite[Theorem~5.1]{DIS}. Let $k\ge2$, let $f$ be the equation \eqref{eq:Cke} of $C=C_{k,\varepsilon}$, and let $D:g=0$ be a line ($\beta=1$) or a smooth conic ($\beta=2$) which is not a component of $C$. Then the curve $C':f'=fg=0$, of degree $d'=d+\beta$, has $\mdr(f')\le1+\beta$. If the line is different from $z=0$, resp.\ the conic does not belong to the family $y^2-2xz-cz^2=0$, then $C'$ satisfies \eqref{eq:DSS} and $\mdr(f')\in\{2,3\}$.

Indeed, recall that $\alpha_C=5/(4k+2\varepsilon)\le\frac58$ and $\alpha_Cd-2=\frac12$, see Proposition~\ref{prop:mdr1} and its proof, and let $p=[1:0:0]$. We first show that $\alpha_{C'}\le\alpha_C$.

Assume first that $g(p)\ne0$. Then the germs $(C',p)$ and $(C,p)$ coincide, hence $\alpha_{C'}\le\alpha_{C',p}=\alpha_{C,p}=\alpha_C$. In fact equality holds, since the new singularities of $C'$ are the intersection points of $D$ with a single component of $C$, of type $A_{2m-1}$ with $m\le4$, whose Arnold exponents $\frac12+\frac1{2m}$ are at least $\frac58$.

Assume now that $g(p)=0$. Then the singularity of $C$ at $p$ changes, and we consider the family of curves $C'_t:f\cdot(tx^{\beta}+g)=0$. For $t\ne0$ small, the curve $tx^\beta+g=0$ does not pass through $p$, so $(C'_t,p)=(C,p)$, while $C'_0=C'$. The semicontinuity of the Arnold exponent \cite{Ku,V8} gives $\alpha_{C',p}\le\alpha_{C,p}$, hence again $\alpha_{C'}\le\alpha_{C',p}\le\alpha_C$.

In both cases
\[
\alpha_{C'}d'-2\le\alpha_Cd-2+\alpha_C\beta\le\frac12+\frac{5\beta}{8}\le\frac74<2 .
\]
If $\mdr(f')\le1$, then \eqref{eq:DSS} holds for $C'$ by Proposition~\ref{prop:mdr1}; otherwise $\mdr(f')\ge2>\alpha_{C'}d'-2$. Under our assumptions on $D$ one has $\mdr(f')\ge2$, hence $\mdr(f')\in\{2,3\}$.

For example, let $C=C_{3,0}$ with $c_i=i$, so $d=6$, $\mdr(f)=1$, $\tau(C)=21$ and $\alpha_C=\frac5{12}$. Adding the line $x+y+z=0$ gives a curve $C'$ of degree $7$ with $\mdr(f')=2$ and $\tau(C')=28$. Its singularities are the point $p$, four nodes and one $A_3$ point, and $\AR(f')$ is generated in degrees $2$ and $4$, so $C'$ is free with exponents $(2,4)$. Here $\alpha_{C'}d'-2=\frac{11}{12}<2$. Adding instead the conic $x^2+y^2+3z^2=0$ gives a curve $C''$ of degree $8$ with $\mdr=3$ and $\tau(C'')=33$. Its singularities are $p$ and $12$ nodes, $\AR$ is minimally generated in degrees $3,6,6$, and $\alpha_{C''}\cdot8-2=\frac43<3$. Both computations were done with \textsc{Singular} \cite{Singular}.

This construction gives many examples, with $\mdr(f')=2$ or $3$, for which \eqref{eq:DSS} holds although some singularity is not weighted homogeneous: when $g(p)\ne0$ and $(k,\varepsilon)\ne(2,0)$, the singularity of $C'$ at $p$ is that of $C_{k,\varepsilon}$.
\end{rk}

\section{Concluding remarks and open questions}\label{sec:final}

We introduced DSS-extremal curves, those for which the lower bound $\mdr(f)\ge\alpha_Cd-2$ is an equality. Curves of pencil type $\Phi(g,h)$ are DSS-extremal, the minimal syzygy being $\nabla g\times\nabla h$. Under transversality their Milnor algebras have an explicit minimal resolution depending only on $(d,\deg g,\deg h)$. For six-cuspidal sextics, extremality separates the torus curves from the general non-torus curves, so Zariski's classical pair is a Ziegler pair detected by $\mdr$.

Several questions remain open.
\begin{enumerate}[(1)]
\item Conjecture~\ref{conj:main}: is every DSS-extremal curve of pencil type? A first test case is Conjecture~\ref{conj:sextic}. Another is the class of curves of degree $6m$ with $6m^2$ cusps, where torus-type curves exist for all $m$ (Corollary~\ref{cor:torus}).
\item For a DSS-extremal curve, is $\AR(f)_{\mdr(f)}$ spanned by the vectors $\nabla g\times\nabla h$ over all pencil structures of $C$? This holds for nodal curves (Proposition~\ref{prop:nodal-basis}), and Example~\ref{ex:H-fails}(2) and the Ceva arrangement (Example~\ref{ex:arr}(1)) support it in general.
\item Which Zariski pairs $(C,C')$ with $C$ of torus type are Ziegler pairs? By Theorem~\ref{thm:A}, this holds whenever the partner $C'$ is not DSS-extremal.
\item Is there a version for hypersurfaces in $\PP^n$, $n\ge3$, using the extension of \eqref{eq:DSS} to isolated weighted homogeneous singularities and the $(n-1)$-fold product $\nabla g_1\times\cdots\times\nabla g_{n-1}$?
\item Conjecture~\ref{conj:nonWH}: does $\mdr(f)\ge\alpha_Cd-2$ hold for every reduced curve? By Lemma~\ref{lem:half} a counterexample has $\mdr(f)<\frac d2-2$. Proposition~\ref{prop:mdr1} and Remark~\ref{rk:add-curve} settle $\mdr(f)\le1$ and give examples with $\mdr(f)=2,3$; the next test case is $\mdr(f)=2$ in general.
\end{enumerate}

\section*{Acknowledgments}
The author is grateful to Alexandru Dimca for many very useful suggestions and discussions.

\end{document}